\documentclass[]{interact}

\usepackage{epstopdf}% To incorporate .eps illustrations using PDFLaTeX, etc.
\usepackage[caption=false]{subfig}% Support for small, `sub' figures and tables
\usepackage{framed}
\usepackage{xcolor}
\usepackage[hidelinks]{hyperref}

\usepackage[numbers,sort&compress]{natbib}% Citation support using natbib.sty
\bibpunct[, ]{[}{]}{,}{n}{,}{,}% Citation support using natbib.sty
\renewcommand\bibfont{\fontsize{10}{12}\selectfont}% Bibliography support using natbib.sty
\makeatletter% @ becomes a letter
\def\NAT@def@citea{\def\@citea{\NAT@separator}}% Suppress spaces between citations using natbib.sty
\makeatother% @ becomes a symbol again

\theoremstyle{plain}% Theorem-like structures provided by amsthm.sty

\theoremstyle{definition}

\theoremstyle{remark}

\begin{document}

% \articletype{ARTICLE TEMPLATE}% Specify the article type or omit as appropriate

\title{Orthogonality relations for conical functions of imaginary order}

\author{
\name{Job Feldbrugge\textsuperscript{a}\thanks{Email: job.feldbrugge@ed.ac.uk} and Nynke M.D.\ Niezink\textsuperscript{b}}
\affil{\textsuperscript{a}Higgs Centre for Theoretical Physics, James Clerk Maxwell Building, Edinburgh EH9 3FD, United Kingdom; \\ \textsuperscript{b}Department of Statistics and Data Science, Carnegie Mellon University,
Baker Hall, Pittsburgh, PA 15213, USA}
}

\maketitle

\begin{abstract}
Orthogonality relations for conical or Mehler functions of imaginary order are derived and expressed in terms of the Dirac delta function. This work extends recently derived orthogonality relations of associated Legendre functions.
\end{abstract}


\section{Introduction}
The associated Legendre (or Ferrers) functions of the first kind $P_\lambda^\mu(x)$ and second kind $Q_\lambda^\mu(x)$, with the complex degree $\lambda$ and complex order $\mu$, are generalizations of the Legendre polynomial $P_l(x)$, where the degree $\lambda = l$ is an integer and the order vanishes, and the associated Legendre polynomial $P_l^m(x)$, where the degree $\lambda=l$ is an integer and the order $\mu=m$ is an integer with $|m| \leq l$. The Legendre polynomial and function are widely used in applied mathematics and theoretical physics. The polynomial, for example, appears in the study of the Newtonian potential \cite{Legendre:1785}, spherical harmonics, and the energy eigenstates of the hydrogen model in quantum mechanics as in, \textit{e.g.}, the textbook \cite{Landau:1981}. Its generalization to complex degrees and orders plays a pivotal role in many applications including the (generalized) Mehler-Fock transform \cite{Nasim:1984}. The Legendre function $P_\lambda^\mu$ with imaginary order $\mu$ arises in the study of modified P\"oschl-Teller potentials in quantum mechanics \cite{Grosche:1998}, backscattering of radiation in plasmas \cite{Kalmykov:2004}, topological black holes \cite{Hutasoit:2009}, two-particle Hamiltonians \cite{Hallnas:2015}, and Yang-Mills matrix models \cite{Karczmarek:2023}.

The associated Legendre polynomial satisfies the orthogonality condition
\begin{align}
    \int_{-1}^{1} \frac{P_l^{m}(x) P_l ^{n}(x)}{1-x^2}dx = \frac{(l+m)!}{m(l-m)!}\delta_{m,n}\,,
\end{align}
for $m,n \neq 0$. Kalmykov and Shvets \cite{Kalmykov:2004} found a similar relation for the associated Legendre function of the first kind with imaginary order, 
\begin{align}
    \int_{-1}^1 \frac{P_1^{iq}(x) P_1^{iq'}(x)}{1-x^2}\mathrm{d}x = \frac{2 \sinh(\pi q)}{q} \delta(q-q')\,,
\end{align}
with real $q,q'$. Hutasiot \textit{et al.}\ \cite{Hutasoit:2009} generalized this relation to integer degree, 
\begin{align}
    \int_{-1}^1 \frac{P_l^{iq}(x) P_l^{iq'}(x)}{1-x^2}\mathrm{d}x = \frac{2 \sinh(\pi q)}{q} \delta(q-q')\,,
\end{align}
for $l=0,1,2,\dots$. Finally, by an elegant algebraic derivation, Bielski \cite{Bielski:2013} demonstrated that these identities are special cases of more general orthogonality relations. In particular, Bielski evaluated the integrals
\begin{align}
    \int_{-1}^1 \frac{P_\lambda^{i q}(x)P_\lambda^{i q'}(x)}{1-x^2}\mathrm{d}x\,, \quad
    \int_{-1}^1 \frac{P_\lambda^{i q}(x)Q_\lambda^{i q'}(x)}{1-x^2}\mathrm{d}x\,,\quad
    \int_{-1}^1 \frac{Q_\lambda^{i q}(x)Q_\lambda^{i q'}(x)}{1-x^2}\mathrm{d}x\,,
    \label{eq:Bielski}
\end{align}
for general complex degree $\lambda$, in terms of the Dirac delta functions $\delta(q-q')$ and $\delta(q+q')$.

In this paper, we evaluate the related orthogonality relations for the conical functions, 
\begin{align}
    I_1^{q,q'} &= \int_{-1}^1 \frac{P_\lambda^{i q}(x)P_\lambda^{i q'}(x)^*}{1-x^2}\mathrm{d}x\,,\qquad  
    J_1^{q,q'} = \int_{-1}^1 \frac{P_\lambda^{i q}(x)P_\lambda^{i q'}(-x)^*}{1-x^2}\mathrm{d}x\,,\label{eq:1}\\
    I_2^{q,q'} &= \int_{-1}^1 \frac{P_\lambda^{i q}(x)Q_\lambda^{i q'}(x)^*}{1-x^2}\mathrm{d}x\,,\qquad
    J_2^{q,q'} = \int_{-1}^1 \frac{P_\lambda^{i q}(x)Q_\lambda^{i q'}(-x)^*}{1-x^2}\mathrm{d}x\,,\label{eq:2}\\
    I_3^{q,q'} &= \int_{-1}^1 \frac{Q_\lambda^{i q}(x)Q_\lambda^{i q'}(x)^*}{1-x^2}\mathrm{d}x\,,\qquad
    J_3^{q,q'} = \int_{-1}^1 \frac{Q_\lambda^{i q}(x)Q_\lambda^{i q'}(-x)^*}{1-x^2}\mathrm{d}x\,,\label{eq:3}
\end{align}
when $\lambda =- \frac{1}{2}+i \nu$ for real $\nu$, in terms of the Dirac delta functions $\delta(q-q')$ and $\delta(q+q')$, following the derivation described in Bielski \cite{Bielski:2013}. We evaluate the $I_i^{q,q'}$ integrals from Bielski's result \eqref{eq:Bielski}, and use these to derive the $J_i^{q,q'}$ integrals.

The Legendre function with degree $\text{Re}[\lambda]=-\frac{1}{2}$ is known as the conical or Mehler function, first analyzed in the study of conics in electrostatics \cite{Mehler:1868}. The conical function plays a prominent role in the Mehler-Fock transformation and, with imaginary degree, forms the energy eigenstates of the modified P\"oschl-Teller or Rosen-Morse barrier in quantum mechanics \cite{Kleinert:1992}. For a detailed study of the conical function with imaginary order, see \cite{Dunster:2013}. 

These orthogonality relations are of special interest when normalizing the continuum spectrum of the modified P\"oschl-Teller model in quantum mechanics. The relations also apply to many other physical models (see, \textit{e.g.}, \cite{Karczmarek:2023}), as the inner product in quantum mechanics includes a complex conjugation, \textit{i.e.},
\begin{align}
    \langle \psi_1|\psi_2 \rangle = \int \psi_1(x)^*\psi_2(x) \mathrm{d}x\,,
\end{align}
using Dirac's bra–ket notation. 

%%%%%%%%%%%%%%%%%%%%%%%%%%
\section{Relevant properties of the associated Legendre function}
To aid the evaluation of the integrals \eqref{eq:1}--\eqref{eq:3}, we briefly summarize several useful properties of associated Legendre functions. The associated Legendre functions $P_\lambda^\mu(x)$ and $Q_\lambda^\mu(x)$ are solutions of the general Legendre equation
\begin{align}
    \frac{\mathrm{d}}{\mathrm{d}x}\left[(1-x^2) \frac{\mathrm{d}w(x)}{\mathrm{d}x}\right] + \left[\lambda(\lambda+1) - \frac{\mu^2}{1-x^2}\right] w(x)= 0\,,
\end{align}
and are related by the equation 
\begin{align}
    Q_\lambda^\mu(x) = \frac{\pi}{2 \sin (\pi \mu)}\left(P_\lambda^\mu(x) \cos(\pi \mu) - \frac{\Gamma(\lambda + \mu + 1)}{\Gamma(\lambda-\mu +1)}P_\lambda^{-\mu}(x)\right)\,,\label{eq:relation1}
\end{align}
where $\Gamma$ denotes the gamma function \cite{Magnus:1967, Abramowitz:1972, Bateman:1995, Gradshteyn:2007, NIST:DLMF} (\textit{e.g.}, equation \href{https://dlmf.nist.gov/14.9}{(14.9.2)} in \cite{NIST:DLMF}). The Legendre function of the first kind is often expressed as 
\begin{align}
    P_\lambda^\mu(x) &= \frac{1}{\Gamma(1-\mu)} \left(\frac{1+x}{1-x}\right)^{\mu/2}\,_2F_1\left(-\lambda,\lambda+1;1-\mu,\frac{1-x}{2}\right)\,,\quad -1<x<1,
\end{align}
in terms of the hypergeometric function 
\begin{align}
    \,_2F_1(a,b;c,x) = \frac{\Gamma(c)}{\Gamma(a)\Gamma(b)}\sum_{n=0}^\infty \frac{\Gamma(a+n)\Gamma(b+n)}{\Gamma(c+n)}\frac{x^n}{n!}\,,
\end{align}
(\textit{e.g.}, equation \href{https://dlmf.nist.gov/14.3}{(14.3.1)} in \cite{NIST:DLMF}). Moreover, the Legendre function satisfies the reflection condition \cite{Magnus:1967, Abramowitz:1972, Bateman:1995, Gradshteyn:2007, NIST:DLMF}
\begin{align}
    P_\lambda^{\mu}(-x) &= P_\lambda^\mu(x) \cos(\pi(\lambda + \mu)) - \frac{2}{\pi}Q_\lambda^\mu(x)\sin(\pi(\lambda+\mu))\,,\label{eq:relation2}
\end{align}
(\textit{e.g.}, equation \href{https://dlmf.nist.gov/14.9}{(14.9.10)} in \cite{NIST:DLMF}). 

For the conical functions, with $\lambda=-\tfrac{1}{2}+i \nu$, conjugation of the degree leaves the Legendre functions invariant,
\begin{align}
    P_{\lambda}^\mu(x) = P_{\lambda^*}^\mu(x)\,,\quad Q_{\lambda}^\mu(x) = Q_{\lambda^*}^\mu(x)\,,
\end{align}
as the Legendre equation is identical, for $\lambda(\lambda+1) = \lambda^*(\lambda^*+1)=-\tfrac{1}{4}-\nu^2$.

% In the asymptotic region $x \to 1^-$, the Legendre function of the first kind approaches
% \begin{align}
%     P_{\lambda}^{\mu}(x) &\overset{x\to 1^-}{\sim} \frac{1}{\Gamma(1-\mu)}\left(\frac{1+x}{1-x}\right)^{\mu/2}.
% \end{align}
% This follows from the expression of the Legendre function in terms of the hypergeometric function. In the limit $x \to -1^+$, the Legendre function of the first kind approaches the asymptotic
% \begin{align}
%     P_{\lambda}^{\mu}(x) \overset{x\to -1^+}{\sim} &- \frac{\sin (\pi \lambda) \Gamma(\mu)}{\pi}\left(\frac{1-x}{1+x}\right)^{\mu/2}\nonumber\\
%     & + \frac{\Gamma(-\mu)}{\Gamma(-\lambda - \mu)\Gamma(1+\lambda - \mu)}\left(\frac{1-x}{1+x}\right)^{-\mu/2}.
% \end{align}

%%%%%%%%%%%%%%%%%%%%%%%%%%
\section{Orthogonality relations of conical functions}
The Legendre function $P_\lambda^{iq}(x)$ and $P_{\lambda}^{iq'}(x)$ satisfy the differential equations
\begin{align}
    &\frac{\mathrm{d}}{\mathrm{d}x}\left[(1-x^2) \frac{\mathrm{d} P_{\lambda}^{iq}(x)}{\mathrm{d}x}\right] + \left[\lambda(\lambda+1) + \frac{q^2}{1-x^2}\right] P_\lambda^{iq}(x)= 0\,,\\
    &\frac{\mathrm{d}}{\mathrm{d}x}\left[(1-x^2) \frac{\mathrm{d} P_{\lambda}^{iq'}(x)}{\mathrm{d}x}\right] + \left[\lambda(\lambda+1) + \frac{q'^2}{1-x^2}\right] P_{\lambda}^{iq'}(x)= 0\,.
\end{align}
Multiplying the first equation by $P_{\lambda}^{iq'}(x)$ and the second equation by $P_\lambda^{iq}(x)$, subtracting the resulting two equations, and integrating the resulting identity from $a$ to $b$ where $-1 < a< b< 1$, we obtain upon integration by parts
\begin{align}
    \int_{a}^b \frac{P_\lambda^{iq}(x)P_{\lambda}^{iq'}(x)}{1-x^2} \mathrm{d}x
    = \frac{\left[P_{\lambda}^{iq}(x) (1-x^2) \frac{\mathrm{d} P_{\lambda}^{iq'}(x)}{\mathrm{d}x}
    -
    P_{\lambda}^{iq'}(x) (1-x^2) \frac{\mathrm{d} P_{\lambda}^{iq}(x)}{\mathrm{d}x}\right]^{b}_{a}}{q^2-q'^2}\,.
\end{align}
Letting $a \to -1$ from above and $b\to 1$ from below, we write the integral of interest as the difference between two limits
\begin{align}
    \int_{-1}^1 \frac{P_\lambda^{iq}(x)P_{\lambda}^{iq'}(x)}{1-x^2} \mathrm{d}x
    &= \lim_{b \to 1^-} \frac{  P_{\lambda}^{iq}(b) (1-b^2) \frac{\mathrm{d} P_{\lambda}^{iq'}(b)}{\mathrm{d}b}- P_{\lambda}^{iq'}(b) (1-b^2) \frac{\mathrm{d} P_{\lambda}^{iq}(b)}{\mathrm{d}b}}{q^2-q'^2}\nonumber\\
    &- \lim_{a \to -1^+} \frac{P_{\lambda}^{iq}(a) (1-a^2) \frac{\mathrm{d} P_{\lambda}^{iq'}(a)}{\mathrm{d}a} - P_{\lambda}^{iq'}(a) (1-a^2) \frac{\mathrm{d} P_{\lambda}^{iq}(a)}{\mathrm{d}a}}{q^2-q'^2}\,.
\end{align}
% The first limit can be evaluated as 
% \begin{align}
%     &\lim_{b \to 1^-} \frac{P_{\lambda}^{iq}(b) (1-b^2) \frac{\mathrm{d} P_{\lambda^*}^{-iq'}(b)}{\mathrm{d}b} - P_{\lambda^*}^{-iq'}(b) (1-b^2) \frac{\mathrm{d} P_{\lambda}^{iq}(b)}{\mathrm{d}b} }{q^2-q'^2}\nonumber \\
%     &=-\frac{i}{\Gamma(1-iq)\Gamma(1+iq')}\lim\limits_{b \to 1^-}\frac{\exp\left[-\frac{i}{2}(q-q')\ln \left(\frac{1-b}{1+b}\right)\right]}{q-q'}\\
%     &=\frac{\pi}{\Gamma(1-iq)\Gamma(1+iq')}\delta(q-q')\\
%     &=\frac{\pi}{\Gamma(1-iq)\Gamma(1+iq)}\delta(q-q')\,,
% \end{align}
% using the identity 
% \begin{align}
%     \lim\limits_{u \to \infty}e^{i u k} &= \lim\limits_{u \to \infty}\cos(u k) + i \lim\limits_{u \to \infty}\sin(u k) =  i \pi k \delta(k)\,,
% \end{align} 
% as the first term vanishes $\lim\limits_{u\to \infty}\cos(k u)=0$ in the distributional sense and the second term
% \begin{align}
%     \lim_{u \to \infty} \sin(u k) = \frac{k}{2} \lim_{u\to \infty} \frac{e^{iuk}- e^{-iuk}}{i k} = \frac{k}{2} \lim_{u \to \infty} \int_{-u}^u e^{ik\alpha}\mathrm{d}\alpha
% \end{align}
% approaches the Fourier representation of the Dirac delta function
% \begin{align}
%     2 \pi \delta(k) = \int_{-\infty}^\infty e^{i k \alpha}\mathrm{d}\alpha\,.
% \end{align}
% The second limit follows analogously, and therefore
%
Bielski \cite{Bielski:2013} evaluates these limits and shows that
\begin{align}
    \int_{a}^b \frac{P_\lambda^{iq}(x)P_{\lambda}^{iq'}(x)}{1-x^2} \mathrm{d}x
    &=- \frac{2 \Gamma(iq)\Gamma(-iq)\sin(\pi \lambda)}{\Gamma(1+\lambda-iq)\Gamma(-\lambda-iq)}\delta(q-q')\nonumber\\
    &= \bigg[\frac{\pi}{\Gamma(1-iq)\Gamma(1+iq)} + \frac{\sin^2(\pi \lambda) \Gamma(iq)\Gamma(-iq)}{\pi} \nonumber\\
    &\phantom{=}+ \frac{\pi \Gamma(iq)\Gamma(-iq)}{\Gamma(1+\lambda - iq)\Gamma(-\lambda - iq)\Gamma(1+\lambda + iq)\Gamma(-\lambda + iq)}\bigg]\,.
\end{align}

Using the fact that for conical functions 
\begin{align}
    P_\lambda^{iq}(x)^*&=P_{\lambda^*}^{-iq}(x)= P_{\lambda}^{-iq}(x)\,,
\end{align}
we find 
\begin{align}
    I_1^{q,q'} 
    &= \frac{\cosh(2 \pi q )+ \cosh(2 \pi \nu)}{q \sinh(\pi q)}\delta(q-q')\nonumber\\
    &\phantom{=} +  \frac{2 \pi \cosh(\pi \nu)}{q \sinh (\pi q) \Gamma(\tfrac{1}{2}-i \nu -iq)\Gamma(\tfrac{1}{2}+ i \nu -iq)}  \delta(q+q')\,.
\end{align}
The integrals $I_2^{q,q'}$ and $I_3^{q,q'}$ follow directly from relation \eqref{eq:relation1}, 
\begin{align}
    I_2^{q,q'} &=\frac{i\pi}{2\sinh(\pi q') }\left[ \cosh( \pi q') I_1^{q,q'}- \frac{\Gamma(\tfrac{1}{2}-i \nu-iq')}{\Gamma(\tfrac{1}{2}-i \nu + iq')}I_1^{q,-q'}\right]\\
    &= \frac{i \pi (\sinh(2 \pi q) + \sinh(2 \pi \nu))}{2 q \sinh (\pi q)}\delta(q-q')
    \nonumber \\
    &\phantom{=}+
    \frac{i \pi^2 \sinh(\pi \nu)}{q \sinh(\pi q) \Gamma(\tfrac{1}{2}-i\nu-iq)\Gamma(\tfrac{1}{2}+i \nu -iq)}\delta(q+q')\,,
\end{align}
and  
\begin{align}
    I_3^{q,q'} &= \frac{i\pi}{2 \sinh( \pi q)} \left[  \frac{\Gamma(\tfrac{1}{2} + i \nu +iq)}{\Gamma(\tfrac{1}{2} + i \nu - iq )}I_2^{-q,q'}-\cosh ( \pi q) I_2^{q,q'}\right]\\
    &= \frac{\pi^2(\cosh(2 \pi q) + \cosh(2\pi \nu))}{4 q \sinh(\pi q)}\delta(q-q') \nonumber\\
    &\phantom{=}+ \frac{\pi^3 \cosh(\pi \nu)}{2q \sinh(\pi q)\Gamma(\tfrac{1}{2} - i \nu-iq)\Gamma(\tfrac{1}{2}  + i \nu- i q)}\delta(q+q')\,.
\end{align}
Using the reflection equation \eqref{eq:relation2}, we obtain the orthogonality relations for $J_1^{q,q'}$,
\begin{align}
    J_1^{q,q'}
    % &= \cos\left(\pi\left(-\tfrac{1}{2} - i \nu - i q'\right)\right) I_1^{q,q'} - \frac{2}{\pi} \sin\left(\pi\left(-\tfrac{1}{2} - i \nu-i q'\right)\right)I_2^{q,q'}\\
    &= -i \sinh(\pi(\nu+q')) I_1^{q,q'} + \frac{2}{\pi} \cosh(\pi(\nu+q'))I_2^{q,q'}\\
    &=\frac{2\pi i}{q \Gamma(\tfrac{1}{2} -i \nu-i q)\Gamma(\tfrac{1}{2} + i \nu- i q )}\delta(q+q')\,.
\end{align}
Again applying equation \eqref{eq:relation1}, we obtain
\begin{align}
    J_2^{q,q'}
    % &= \frac{\pi}{2 \sin\left(\pi\left(-\tfrac{1}{2}-i\nu-iq'\right)\right)}\left[\cos\left(\pi\left(-\tfrac{1}{2}-i \nu - iq'\right)\right)J_1^{q,q'} - I_1^{q,q'}\right]\\
    % &= \frac{\pi}{2 \cosh(\pi(\nu+q'))}\left[i\sinh(\pi(\nu+q'))J_1^{q,q'} + I_1^{q,q'}\right]\\
    &=\frac{i \pi}{2 \sinh(\pi q')}\left[\cosh(\pi q')J_1^{q,q'}-\frac{\Gamma(\tfrac{1}{2} - i \nu - i q')}{\Gamma(\tfrac{1}{2} - i \nu + i q')}J_1^{q,-q'}\right]\\
    &=\frac{\cosh(\pi (q-\nu))}{q \sinh(\pi q)}\delta(q-q') \nonumber\\
    &\phantom{=}+ \frac{\pi^2}{q \tanh(\pi q)\Gamma(\tfrac{1}{2} -i \nu- i q )\Gamma(\tfrac{1}{2}  + i \nu- i q)}\delta(q+q')\,,
\end{align}
and 
\begin{align}
    J_3^{q,q'}&= \frac{i\pi}{2 \sinh(\pi  q)}\left[\frac{\Gamma(\tfrac{1}{2}+i \nu +iq)}{\Gamma(\tfrac{1}{2}+i \nu -iq)}J_2^{-q,q'} - \cosh (\pi q)J_2^{q,q'}\right]\\
    &= -\frac{i \pi^3 }{2 q \Gamma(\tfrac{1}{2}- i \nu-i q )\Gamma(\tfrac{1}{2}+ i \nu - i q )} \delta(q + q')\,.
\end{align}
Note that Bielski's formulas do not directly imply the $J_i^{q,q'}$ integrals.

These expressions satisfy the following conjugation conditions which we can obtain from their definitions,
\begin{align}
    &\left(I_1^{q,q'}\right)^* = I_1^{q',q}\,,  &&\left(I_3^{q,q'}\right)^*=I_3^{q',q}\,,\\
    &\left(I_2^{q,q'}\right)^* = I_2^{-q,-q'}\,,  &&\left( J_2^{q,q'}\right)^* =J_2^{-q,-q'}\,,\\
    &\left(J_1^{q,q'}\right)^* = J_1^{q',q} \,,  &&\left(J_3^{q,q'}\right)^*=J_3^{q',q}\,.
\end{align}
From the derived equations, we find that $ \left(I_2^{q,q'}\right)^*=  - I_2^{q',q}$ and $ \left(J_2^{q,q'}\right)^*=  J_2^{q',q}$.

\section*{Acknowledgements}
The work of JF is supported by the STFC Consolidated Grant ‘Particle Physics at the Higgs Centre,’ and, respectively, by a Higgs Fellowship and the Higgs Chair of Theoretical Physics at the University of Edinburgh. 

For the purpose of open access, the author has applied a Creative Commons Attribution (CC BY) license to any Author Accepted Manuscript version arising from this submission.

\bibliographystyle{tfnlm}
% \bibliography{interactnlmsample}
\bibliography{Library}

\end{document}